\documentclass{amsart}
\usepackage{amssymb}
\usepackage{amsmath}
\usepackage{amsfonts}

\newtheorem{theorem}{Theorem}
\theoremstyle{plain}

\newtheorem{corollary}{Corollary}

\newtheorem{lemma}{Lemma}

\newtheorem{proposition}{Proposition}
\newtheorem{remark}{Remark}

\numberwithin{equation}{section}
\input{tcilatex}

\begin{document}
\title[Properties of modified Riemannian extensions]{Properties of modified
Riemannian extensions}
\author{Aydin GEZER}
\address{Ataturk University, Faculty of Science, Department of Mathematics,
25240, Erzurum-Turkey.}
\email{agezer@atauni.edu.tr}
\author{Lokman BILEN}
\address{Igdir University, Igdir Vocational School, 76000, Igdir-Turkey.}
\email{lokman.bilen@igdir.edu.tr}
\author{Ali CAKMAK}
\address{Ataturk University, Faculty of Science, Department of Mathematics,
25240, Erzurum-Turkey.}
\email{ali.cakmak@atauni.edu.tr}

\begin{abstract}
Let $M$ be an $n-$dimensional differentiable manifold with a symmetric
connection $\nabla $ and $T^{\ast }M$ be its cotangent bundle. In this
paper, we study some properties of the modified Riemannian extension $%
\widetilde{g}_{\nabla ,c}$ on $T^{\ast }M$ defined by means of a symmetric $%
(0,2)$-tensor field $c$ on $M.$ We get the conditions under which $T^{\ast }M
$ endowed with the horizontal lift $^{H}J$ of an almost complex structure $J$
and with the metric $\widetilde{g}_{\nabla ,c}$ is a K\"{a}hler-Norden
manifold. Also curvature properties of the Levi-Civita connection and
another metric connection of the metric $\widetilde{g}_{\nabla ,c}$ are
presented.
\end{abstract}

\maketitle

\section{\protect\bigskip \textbf{Introduction}}

\noindent Let $M$ be an $n-$dimensional differentiable manifold and $T^{\ast
}M$ be its cotangent bundle. There is a well-known natural construction
which yields, for any affine connection $\nabla $ on $M$, a
pseudo-Riemannian metric $\widetilde{g}_{\nabla }$ on $T^{\ast }M.$ The
metric $\widetilde{g}_{\nabla }$ is called Riemannian extension of $\nabla .$
Riemannian extensions were originally defined by Patterson and Walker \cite%
{Patterson} and further investigated in Afifi \cite{Afifi}, thus relating
pseudo-Riemannian properties of $T^{\ast }M$ with the affine structure of
the base manifold $(M,\nabla )$. Moreover, Riemannian extensions were also
considered in Garcia-Rio et al. \cite{Garcia1} in relation to Osserman
manifolds (see also Derdzinski \cite{Derdzinski}). Riemannian extensions
provide a link between affine and pseudo-Riemannian geometries. Some
properties of the affine connection $\nabla $ can be investigated by means
of the corresponding properties of the Riemannian extension $\widetilde{g}%
_{\nabla }$. For instance, $\nabla $ is projectively flat if and only if $%
\widetilde{g}_{\nabla }$ is locally conformally flat \cite{Afifi}. For
Riemannian extensions, also see \cite%
{Aslanci,Dryuma,Honda,Kowalski,Mok,Sekizawa1,Toomanian,Vanhecke,Willmore}.
In \cite{Calvino1,Calvino2}, the authors introduced a modification of the
usual Riemannian extensions which is called modified Riemannian extension.

Almost complex Norden manifolds are among the most important geometrical
structures which can be considered on a manifold. Let $M_{2k}$ be a $2k$%
-dimensional differentiable manifold endowed with an almost complex
structure $J$ and a pseudo-Riemannian metric $g$ of signature $(k,k)$ such
that $g(JX,Y)=g(X,JY)$, i.e. $g$ is pure with respect to $J$ for arbitrary
vector fields $X$ and $Y$ on $M_{2k}$. Then the metric $g$ is called Norden
metric. Norden metrics are referred to as anti-Hermitian metrics or $B$%
-metrics. The study of such manifolds is interesting because there exists a
difference between the geometry of a $2k-$dimensional almost complex
manifold with Hermitian metric and the geometry of a $2k-$dimensional almost
complex manifold with Norden metric. A notable difference between Norden
metrics and Hermitian metrics is that $G(X,Y)=g(X,JY)$ is another Norden
metric, rather than a differential $2$-form. Some authors considered almost
complex Norden structures on the cotangent bundle \cite%
{Druta1,Oproiu1,Oproiu2}.

In this paper we will use a deformation of the Riemannian extension on the
cotangent bundle $T^{\ast }M$ over $(M,\nabla )$ by means of a symmetric
tensor field $c$ on $M$, where $\nabla $ is a symmetric affine connection on 
$M$. The metric is so-called modified Riemannian extenson. In section 3, in
the particular case where $\nabla $ is the Levi-Civita connection on a
Riemannian manifold $(M,g)$, we get the conditions under which the almost
complex manifold with Norden metric $(T^{\ast }M,^{H}J,\widetilde{g}_{\nabla
,c})$ is a K\"{a}hler-Norden manifold, where $^{H}J$ is the horizontal lift
of an almost complex structure $J$ and $\widetilde{g}_{\nabla ,c}$ is the
modified Riemannian extension. In section 4 and 5, we show that the
geometric properties of the Levi-Civita connection and another metric
connection of the modified Riemannian extension $\widetilde{g}_{\nabla ,c}$.

Throughout this paper, all manifolds, tensor fields and connections are
always assumed to be differentiable of class $C^{\infty }$. Also, we denote
by $\Im _{q}^{p}(M)$ the set of all tensor fields of type $(p,q)$ on $M$,
and by $\Im _{q}^{p}(T^{\ast }M)$ the corresponding set on the cotangent
bundle $T^{\ast }M$. The Einstein summation convention is used, the range of
the indices $i,j,s$ being always $\{1,2,...,n\}$

\section{Preliminaries}

\subsection{The cotangent bundle}

Let $M$ be an $n-$dimensional smooth manifold and denote by $\pi :T^{\ast
}M\rightarrow M$ its cotangent bundle with fibres the cotangent spaces to $%
M. $ Then ${}T^{\ast }M$ is a $2n-$dimensional smooth manifold and some
local charts induced naturally from local charts on $M,$ may be used.
Namely, a system of local coordinates $\left( U,x^{i}\right) ,\mathrm{\;}%
i=1,...,n$ in $M$ induces on ${}T^{\ast }M$ a system of local coordinates $%
\left( \pi ^{-1}\left( U\right) ,\mathrm{\;}x^{i},\mathrm{\;}x^{\overline{i}%
}=p_{i}\right) ,\mathrm{\;}\overline{i}=n+i=n+1,...,2n$, where $x^{\overline{%
i}}=p_{i}$ is the components of covectors $p$ in each cotangent space $%
{}T_{x}^{\ast }M,\mathrm{\;}x\in U$ with respect to the natural coframe $%
\left\{ dx^{i}\right\} $.

Let $X=X^{i}\frac{\partial }{\partial x^{i}}$ and $\omega =\omega _{i}dx^{i}$
be the local expressions in $U$ of a vector field $X$ \ and a covector
(1-form) field $\omega $ on $M$, respectively. Then the vertical lift $%
^{V}\omega $ of $\omega $, the horizontal lift $^{H}X$ and the complete lift 
$^{C}X$ of $X$ are given, with respect to the induced coordinates,
by\noindent 
\begin{equation}
^{V}\omega =\omega _{i}\partial _{\overline{i}},  \label{AB2.1}
\end{equation}

\begin{equation}
^{H}X=X^{i}\partial _{i}+p_{h}\Gamma _{ij}^{h}X^{j}\partial _{\overline{i}}
\label{AB2.2}
\end{equation}%
and%
\begin{equation*}
^{C}X=X^{i}\partial _{i}-p_{h}\partial _{i}X^{h}\partial _{\overline{i}},
\end{equation*}%
where $\partial _{i}=\frac{\partial }{\partial x^{i}}$, $\partial _{%
\overline{i}}=\frac{\partial }{\partial x^{\overline{i}}}$ and $\Gamma
_{ij}^{h}$ are the coefficients of a symmetric (torsion-free) affine
connection $\nabla $ in $M.$

The Lie bracket operation of vertical and horizontal vector fields on $%
T^{\ast }M$ is given by the formulas%
\begin{equation}
\left\{ 
\begin{array}{l}
{\left[ ^{H}X,^{H}Y\right] ={}^{H}\left[ X,Y\right] +^{V}\left( p\circ
R(X,Y)\right) } \\ 
{\left[ ^{H}X,^{V}\omega \right] ={}^{V}\left( \nabla _{X}\omega \right) }
\\ 
{\left[ ^{V}\theta ,^{V}\omega \right] =0}%
\end{array}%
\right.  \label{AB2.3}
\end{equation}%
for any $X,$ $Y$ $\Im _{0}^{1}(M)$ and $\theta $, $\omega \in \Im
_{1}^{0}(M) $, where $R$ is the curvature tensor of the symmetric connection 
$\nabla $ defined by $R\left( X,Y\right) =\left[ \nabla _{X},\nabla _{Y}%
\right] -\nabla _{\left[ X,Y\right] }$ (for details, see \cite%
{YanoIshihara:DiffGeo})$.$

\subsection{Expressions in the adapted frame}

We insert the adapted frame which allows the tensor calculus to be
efficiently done in $T^{\ast }M.$ With the symmetric affine connection $%
\nabla $ in $M$, we can introduce adapted frames on each induced coordinate
neighborhood $\pi ^{-1}(U)$ of $T^{\ast }M$. In each local chart $U\subset M$%
, we write $X_{(j)}=\dfrac{\partial }{\partial x^{j}},$ $\theta
^{(j)}=dx^{j},$ $j=1,...,n.$ Then from (\ref{AB2.1}) and (\ref{AB2.2}), we
see that these vector fields have, respectively, local expressions 
\begin{equation*}
^{H}X_{(j)}=\partial _{j}+p_{a}\Gamma _{hj}^{a}\partial _{\overline{h}}
\end{equation*}%
\begin{equation*}
^{V}\theta ^{(j)}=\partial _{\overline{j}}
\end{equation*}%
with respect to the natural frame $\left\{ \partial _{j},\partial _{%
\overline{j}}\right\} $. These $2n$ vector fields are linearly independent
and they generate the horizontal distribution of $\nabla $ and the vertical
distribution of $T^{\ast }M$, respectively. The set $\left\{
^{H}X_{(j)},^{V}\theta ^{(j)}\right\} $ is called the frame adapted to the
connection $\nabla $ in $\pi ^{-1}(U)\subset T^{\ast }M$. By denoting%
\begin{eqnarray}
E_{j} &=&^{H}X_{(j)},  \label{AB2.4} \\
E_{\overline{j}} &=&^{V}\theta ^{(j)},  \notag
\end{eqnarray}%
we can write the adapted frame as $\left\{ E_{\alpha }\right\} =\left\{
E_{j},E_{\overline{j}}\right\} $. The indices $\alpha ,\beta ,\gamma
,...=1,...,2n$ indicate the indices with respect to the adapted frame.

Using (\ref{AB2.1}), (\ref{AB2.2}) and (\ref{AB2.4}), we have%
\begin{equation}
^{V}\omega =\left( 
\begin{array}{l}
0 \\ 
\omega _{j}%
\end{array}%
\right) ,  \label{AB2.5}
\end{equation}%
and

\begin{equation}
^{H}X=\left( 
\begin{array}{l}
X^{j} \\ 
0%
\end{array}%
\right)  \label{AB2.6}
\end{equation}%
with respect to the adapted frame $\left\{ E_{\alpha }\right\} $ (for
details, see \cite{YanoIshihara:DiffGeo}). By the straightforward
calculations, we have the lemma below.

\begin{lemma}
\label{Lemma1}The Lie brackets of the adapted frame of $T^{\ast }M$ satisfy
the following identities:%
\begin{eqnarray*}
\left[ E_{i},E_{j}\right]  &=&p_{s}R_{ijl}^{\text{ \ \ }s}E_{\overline{l}},
\\
\left[ E_{i},E_{\overline{j}}\right]  &=&-\Gamma _{il}^{j}E_{\overline{l}},
\\
\left[ E_{\overline{i}},E_{\overline{j}}\right]  &=&0,
\end{eqnarray*}%
where $R_{ijl}^{\text{ \ \ }s}$ denote the components of the curvature
tensor of the symmetric connection $\nabla $ on $M$.\bigskip 
\end{lemma}

\section{\protect\bigskip K\"{a}hler-Norden structures on the cotangent
bundle}

An almost complex Norden manifold $(M,J,g)$ (an almost complex manifold with
a Norden metric) is defined to be a real $2n$-dimensional differentiable
manifold $M$, i.e. $J$ is an almost complex structure and $g$ is a
pseudo-Riemannian metric of neutral signature $(n,n)$ on $M$ such that:%
\begin{equation*}
J^{2}X=-X,\text{ }g(JX,Y)=g(X,JY)
\end{equation*}%
for all $X,Y\in \Im _{0}^{1}(M).$ An K\"{a}hler-Norden (anti-K\"{a}hler)
manifold can be defined as a triple $(M,J,g)$ which consists of a smooth
manifold $M$ endowed with an almost complex structure $J$ and a Norden
metric $g$ such that $\nabla J=0$, where $\nabla $ is the Levi-Civita
connection of $g$. It is well known that the condition $\nabla J=0$ is
equivalent to $\mathrm{C}$-holomorphicity (analyticity) of the Norden metric 
$g$ \cite{Iscan}, i.e. $\Phi _{J}g=0$, where $\Phi _{J}$ is the Tachibana
operator \cite{YanoAko:1968,Tachibana}: $(\Phi
_{J}g)(X,Y,Z)=(JX)({}g(Y,Z))-X(g(JY,Z))+g((L_{Y}J)X,Z)+g(Y,(L_{Z}$ $J)X)\,.$
Also note that $G(Y,Z)=(g\circ J)(Y,Z)=g(JY,Z)$ is the twin Norden metric.
Since in dimension $2$ an K\"{a}hler-Norden manifold is flat, we assume in
the sequel that $\dim \,M\geq 4$.

Next, for a given symmetric connection $\nabla $ on $M$, the cotangent
bundle $T^{\ast }M$ may be equipped with a pseudo-Riemannian metric $%
\widetilde{g}_{\nabla }$ of signature $(n,n)$: the Riemannian extension of $%
\nabla $ \cite{Patterson}, given by 
\begin{equation*}
\widetilde{g}_{\nabla }(^{C}X,^{C}Y)=-\gamma (\nabla _{X}Y+\nabla _{Y}X)
\end{equation*}%
where $^{C}X,^{C}Y$ denote the complete lifts to $T^{\ast }M$ of vector
fields $X,Y$ on $M$. Moreover for any vector field $Z$ on $M$, $%
Z=Z^{i}\partial _{i}$, $\gamma Z$ is the function on $T^{\ast }M$ defined by 
$\gamma Z=p_{i}Z^{i}$. The Riemannian extension is expressed by%
\begin{equation*}
\widetilde{g}_{\nabla }=\left( 
\begin{array}{cc}
-2p_{h}\Gamma _{ij}^{h} & \delta _{j}^{i} \\ 
\delta _{i}^{j} & 0%
\end{array}%
\right) 
\end{equation*}%
with respect to the natural frame.

Now we give a deformation of the Riemannian extension above by means of a
symmetric $(0,2)$-tensor field $c$ on $M$, which the metric is called
modified Riemannian extension. The modified Riemannian extension is
expressed by%
\begin{equation}
{}\widetilde{g}_{\nabla ,c}=g_{\nabla }+\pi ^{\ast }c=\left( 
\begin{array}{cc}
-2p_{h}\Gamma _{ij}^{h}+c_{ij} & \delta _{j}^{i} \\ 
\delta _{i}^{j} & 0%
\end{array}%
\right) .  \label{AB3.1}
\end{equation}%
with respect to the natural frame. It follows that the signature of ${}%
\widetilde{g}_{\nabla ,c}$ is $(n,n).$

Denote by $\nabla $ the Levi-Civita connection of a pseudo-Riemannian metric 
$g$. In the section, we will consider $T^{\ast }M$ equipped with the metric $%
{}\widetilde{g}_{\nabla ,c}$ over a pseudo-Riemannian manifold $(M,g).$
Since the vector fields $^{H}X$ and $^{V}\omega $ span the module of vector
fields on $T^{\ast }M$, any tensor field is determined on $T^{\ast }M$ by
their actions on $^{H}X$ and $^{V}\omega $. The modified Riemannian
extension ${}\widetilde{g}_{\nabla ,c}$ has the following properties%
\begin{eqnarray}
{}\widetilde{g}_{\nabla ,c}(^{H}X,^{H}Y) &=&c(X,Y)  \label{AB3.2} \\
{}\widetilde{g}_{\nabla ,c}(^{H}X,^{V}\omega ) &=&g_{\nabla ,c}(^{V}\omega
,^{H}X)=\omega (X)  \notag \\
{}\widetilde{g}_{\nabla ,c}(^{V}\omega ,^{V}\theta ) &=&0  \notag
\end{eqnarray}%
for all $X,Y\in \Im _{0}^{1}(M)$ and $\omega ,\theta \in \Im _{1}^{0}(M),$
which characterize $\widetilde{g}_{\nabla ,c}.$

The horizontal lift of a $(1,1)$-tensor field $J$ to $T^{\ast }M$ is defined
by:%
\begin{eqnarray}
^{H}J(^{H}X) &=&^{H}(JX)  \label{AB3.3} \\
^{H}J({}^{V}\omega ) &=&^{V}(\omega \circ J)  \notag
\end{eqnarray}%
for any $X\in \Im _{0}^{1}\left( M\right) $ and $\omega \in \Im
_{1}^{0}\left( M\right) $. Moreover, it is well known that \i f $J$ is an
almost complex structure on $(M,g)$, then its horizontal lift $^{H}J$ is an
almost complex structure on $T^{\ast }M$ \cite{YanoIshihara:DiffGeo}.

Putting

\begin{equation*}
A\left( \widetilde{X},\widetilde{Y}\right) ={}{}\widetilde{g}_{\nabla
,c}\left( ^{H}J\widetilde{X},\widetilde{Y}\right) -{}\widetilde{g}_{\nabla
,c}\left( \widetilde{X},^{H}J\widetilde{Y}\right)
\end{equation*}%
for any $\widetilde{X},\widetilde{Y}\in \Im _{0}^{1}\left( T^{\ast }M\right) 
$. For all vector fields $\widetilde{X}$ and $\widetilde{Y}$ which are of
the form ${}^{V}\omega ,\,{}^{V}\theta $ or ${}^{H}X,\,{}^{H}Y$, from (\ref%
{AB3.2}) and (\ref{AB3.3}), we have 
\begin{eqnarray*}
A\left( ^{H}X,^{H}Y\right) &=&{}\widetilde{g}_{\nabla ,c}\left(
^{H}J(^{H}X),^{H}Y\right) -{}{}\widetilde{g}_{\nabla ,c}\left(
^{H}X,^{H}J(^{H}Y)\right) \\
&=&{}{}\widetilde{g}_{\nabla ,c}\left( ^{H}(JX),^{H}Y\right) -{}{}\widetilde{%
g}_{\nabla ,c}\left( ^{H}X,^{H}(JY)\right) \\
&=&c(JX,Y)-c(X,JY) \\
A\left( ^{H}X,^{V}\theta \right) &=&{}{}\widetilde{g}_{\nabla ,c}\left(
^{H}J(^{H}X),^{V}\theta \right) -{}{}\widetilde{g}_{\nabla ,c}\left(
^{H}X,^{H}J(^{V}\theta )\right) \\
{} &=&{}\widetilde{g}_{\nabla ,c}\left( ^{H}(JX),^{V}\theta \right) -{}%
\widetilde{g}_{\nabla ,c}\left( ^{H}X,^{V}(\theta \circ J)\right) \\
&=&\theta (JX)-(\theta \circ J)(X) \\
A\left( ^{V}\omega ,^{H}Y\right) &=&{}\widetilde{g}_{\nabla ,c}\left(
^{H}J(^{V}\omega ),^{H}Y\right) -{}\widetilde{g}_{\nabla ,c}\left(
^{V}\omega ,^{H}J(^{H}X)\right) \\
&=&(\omega \circ J)(Y)-\omega (JX) \\
A\left( ^{V}\omega ,^{V}\theta \right) &=&{}\widetilde{g}_{\nabla ,c}\left(
^{H}J(^{V}\omega ),^{V}\theta \right) -{}\widetilde{g}_{\nabla ,c}\left(
^{V}\omega ,^{H}J(^{V}\theta )\right) \\
&=&{}\widetilde{g}_{\nabla ,c}\left( ^{H}J(^{V}\omega ),^{V}\theta \right)
-{}\widetilde{g}_{\nabla ,c}\left( ^{V}\omega ,^{H}J(^{V}\theta )\right) \\
&=&{}\widetilde{g}_{\nabla ,c}\left( ^{V}(\omega \circ J),^{V}\theta \right)
-{}{}\widetilde{g}_{\nabla ,c}\left( ^{V}\omega ,^{V}(\theta \circ J)\right)
\\
&=&0.
\end{eqnarray*}%
If $g$ and $c$ is pure with respect to $J$, we say that $A\left( \tilde{X},%
\tilde{Y}\right) =0$, i.e. ${}g_{\nabla ,c}$ is pure with respect to $^{H}J.$

We now are interested in the holomorphy property of the metric $g_{\nabla
,c} $ with respect to $^{H}$ $J$. We calculate 
\begin{eqnarray*}
(\Phi _{^{H}J}{}\widetilde{g}_{\nabla ,c})(\tilde{X},\tilde{Y},\tilde{Z})
&=&(^{H}J\tilde{X})({}\widetilde{g}_{\nabla ,c}(\tilde{Y},\tilde{Z}))-\tilde{%
X}(\widetilde{g}_{\nabla ,c}(^{H}J\tilde{Y},\tilde{Z})) \\
&+&{}\widetilde{g}_{\nabla ,c}((L_{\tilde{Y}}\text{ }^{H}J)\tilde{X},\tilde{Z%
})+\widetilde{g}_{\nabla ,c}(\tilde{Y},(L_{\tilde{Z}}\text{ }^{H}J)\tilde{X})
\end{eqnarray*}%
for all $\tilde{X},\tilde{Y},\tilde{Z}\in \Im _{0}^{1}(T^{\ast }M)$. Then we
obtain the following equations 
\begin{eqnarray*}
(\Phi _{^{H}J}{}\widetilde{g}_{\nabla ,c})({}^{V}\omega ,{}^{V}\theta
,{}^{H}Z) &=&0, \\
(\Phi _{^{H}J}{}\widetilde{g}_{\nabla ,c})({}^{V}\omega ,{}^{V}\theta
,{}^{V}\sigma ) &=&0, \\
(\Phi _{^{H}J}{}\widetilde{g}_{\nabla ,c})({}^{V}\omega
,{}^{H}Y,{}^{V}\sigma ) &=&{}0, \\
(\Phi _{^{H}J}{}{}\widetilde{g}_{\nabla ,c})({}^{V}\omega ,{}^{H}Y,{}^{H}Z)
&=&(\omega \circ \nabla _{Y}J)(Z)+(\omega \circ \nabla _{Z}J)(Y), \\
(\Phi _{^{H}J}{}\widetilde{g}_{\nabla ,c})({}^{H}X,{}^{V}\omega ,{}^{H}Z)
&=&(\Phi _{J}{}g)({}X,{}\widetilde{\omega },{}Z)-g((\nabla _{\widetilde{%
\omega }}J)X,Z), \\
(\Phi _{^{H}J}{}\widetilde{g}_{\nabla ,c})({}^{H}X,{}^{V}\omega
,{}^{V}\sigma ) &=&0, \\
(\Phi _{^{H}J}{}\widetilde{g}_{\nabla ,c})({}^{H}X,{}^{H}Y,{}^{H}Z) &=&(\Phi
_{J}{}c)({}X,{}Y,{}Z)) \\
&&+(p\circ R(Y,JX)-p\circ R(Y,X)J)(Z) \\
&&+(p\circ R(Z,JX)-p\circ R(Z,X)J)(Y), \\
(\Phi _{^{H}J}{}\widetilde{g}_{\nabla ,c})({}^{H}X,{}^{H}Y,{}^{V}\sigma )
&=&(\Phi _{J}{}g)({}X,{}Y,{}\widetilde{\sigma })-g(Y,(\nabla _{\widetilde{%
\sigma }}J)X),
\end{eqnarray*}%
where $\widetilde{\omega }=g^{-1}\circ \omega =g^{ij}\omega _{j}$ is the
associated vector field of $\omega $. On the other hand, we know that the
Riemannian curvature $R$ of K\"{a}hler-Norden manifolds is totally pure.
Also, the condition $\Phi _{J}g=0$ is equivalent to $\nabla J=0$, where $%
\nabla $ is the Levi-Civita connection of $g$. Hence, we say following
result.

\begin{theorem}
Let $(M,J,g)$ is a K\"{a}hler-Norden manifold. Then $T^{\ast }M$ is a K\"{a}%
hler-Norden manifold equipped with the metric ${}\widetilde{g}_{\nabla ,c}$
and the almost complex structure $^{H}J$ if and only if the$\ $symmetric $%
(0,2)-$tensor field $c$ on $M$ is a holomorphic tensor field with respect to
the almost complex structure $J$.
\end{theorem}

\section{Curvature properties of the Levi-Civita connection of the modified
Riemannian extension $\protect\widetilde{g}_{\protect\nabla ,c}$}

From now on we will consider $T^{\ast }M$ equipped with the modified
Riemannian extension $\widetilde{g}_{\nabla ,c}$ for a given symmetric
connection $\nabla $ on $M.$ By virtue of (\ref{AB2.5}) and (\ref{AB2.6}),
the modified Riemannian extension $\widetilde{g}_{\nabla ,c}$ and its
inverse $\overline{g}_{\nabla ,c}$ have the following components with
respect to the adapted frame $\left\{ E_{\alpha }\right\} $:%
\begin{equation}
(\widetilde{g}_{\nabla ,c})_{\beta \gamma }=\left( 
\begin{array}{cc}
{c}_{ij} & {\delta }_{i}^{j} \\ 
{\delta }_{j}^{i} & {0}%
\end{array}%
\right) .\text{ }  \label{AB4.1}
\end{equation}%
\begin{equation}
(\overline{g}_{\nabla ,c})^{\beta \gamma }=\left( 
\begin{array}{cc}
{0} & {\delta }_{j}^{i} \\ 
{\delta }_{i}^{j} & -{c}_{ij}%
\end{array}%
\right) .  \label{AB4.2}
\end{equation}

The Levi-Civita connection \bigskip ${}\widetilde{{\nabla }}$ of $\widetilde{%
g}_{\nabla ,c}$ is characterized by the Koszul formula:%
\begin{eqnarray*}
2\widetilde{g}_{\nabla ,c}({}\widetilde{{\nabla }}_{\widetilde{X}}\widetilde{%
Y},\widetilde{Z}) &=&\widetilde{X}(\widetilde{g}_{\nabla ,c}(\widetilde{Y},%
\widetilde{Z}))+\widetilde{Y}(\widetilde{g}_{\nabla ,c}(\widetilde{Z},%
\widetilde{X}))-\widetilde{Z}(\widetilde{g}_{\nabla ,c}(\widetilde{X},%
\widetilde{Y})) \\
-\widetilde{g}_{\nabla ,c}(\widetilde{X},[\widetilde{Y},\widetilde{Z}]) &+&%
\widetilde{g}_{\nabla ,c}(\widetilde{Y},[\widetilde{Z},\widetilde{X}])+\text{
}\widetilde{g}_{\nabla ,c}(\widetilde{Z},[\widetilde{X},\widetilde{Y}])
\end{eqnarray*}%
for all vector fields $\widetilde{X},\widetilde{Y}$ and $\widetilde{Z}$ on $%
T^{\ast }M$. One can verify the Koszul formula for pairs $\widetilde{X}=$ $%
E_{i},E_{\overline{i}}$ and $\widetilde{Y}=$ $E_{j},E_{\overline{j}}$ and $%
\widetilde{Z}=$ $E_{k},E_{\overline{k}}$. By using (\ref{AB4.1}), Lemma \ref%
{Lemma1}, we obtain the following result.

\begin{proposition}
\label{propo1}The Levi-Civita connection \bigskip ${}\widetilde{{\nabla }}$
of the modified Riemannian extension $\widetilde{g}_{\nabla ,c}$ is given by%
\begin{eqnarray*}
\widetilde{\nabla }_{E_{\overline{i}}}E_{\overline{j}} &=&0,\text{ }%
\widetilde{\nabla }_{E_{\overline{i}}}E_{j}=0, \\
\widetilde{\nabla }_{E_{i}}E_{\overline{j}} &=&-\Gamma _{ih}^{j}E_{\overline{%
h}}, \\
\widetilde{\nabla }_{E_{i}}E_{j} &=&\Gamma _{ij}^{h}E_{h}+\{p_{s}R_{hji}^{%
\text{ \ \ \ \ }s}+\frac{1}{2}(\nabla _{i}c_{jh}+\nabla _{j}c_{ih}-\nabla
_{h}c_{ij})\}E_{\overline{h}}
\end{eqnarray*}%
where $R_{hji}^{\text{ \ \ \ \ }s}$ are the local coordinate components of
the curvature tensor field $R$ of the symmetric connection $\nabla $ on $M.$
\end{proposition}

An important geometric problem is to find the geodesics on the smooth
manifolds with respect to the Riemannian metrics. Let $C$ be a curve in $M$
expressed locally by $x^{h}=x^{h}(t)$. We define a curve $\widetilde{C}$ in $%
T^{\ast }M$ by%
\begin{equation}
\left\{ 
\begin{array}{l}
x^{h}=x^{h}(t) \\ 
x^{\overline{h}}\overset{def}{=}p_{h}=\omega _{h}(t)%
\end{array}%
\right. \hspace{5mm}  \label{G.0}
\end{equation}%
where $\omega _{h}(t)$ is a covector field along $C.$ The geodesics of the
connection $\widetilde{\nabla }$ is given by the differential equations%
\begin{equation}
\dfrac{\delta ^{2}x^{A}}{dt^{2}}=\dfrac{d^{2}x^{A}}{dt^{2}}+\text{ }%
\widetilde{\Gamma }_{CB}^{A}\dfrac{dx^{C}}{dt}\dfrac{dx^{B}}{dt}=0
\label{G.1}
\end{equation}%
with respect to the induced coordinates $(x^{h},x^{\overline{h}})$, where $t$
is the arc length of a curve in $T^{\ast }M$.

We write down the form equivalent to (\ref{G.1}), namely,%
\begin{equation*}
\dfrac{d}{dt}(\dfrac{\theta ^{\alpha }}{dt})+\widetilde{\Gamma }_{\gamma
\beta }^{\alpha }\dfrac{\theta ^{\gamma }}{dt}\dfrac{\theta ^{\beta }}{dt}=0
\end{equation*}%
with respect to adapted frame $\left\{ E_{\alpha }\right\} $, where%
\begin{eqnarray*}
\dfrac{\theta ^{h}}{dt} &=&\dfrac{dx^{h}}{dt}, \\
\dfrac{\theta ^{\bar{h}}}{dt} &=&\dfrac{\delta p_{h}}{dt}
\end{eqnarray*}%
along a curve $x^{A}=x^{A}(t)$ in $T^{\ast }M$ \cite{YanoIshihara:DiffGeo}.
Taking account of Proposition \ref{propo1}, then we have%
\begin{equation}
\left\{ 
\begin{array}{l}
(a)\mathit{\ \ \ \ }\dfrac{d^{2}x^{h}}{dt^{2}}+\Gamma _{ij}^{h}\dfrac{dx^{i}%
}{dt}\dfrac{dx^{j}}{dt}=0, \\ 
(b)\mathit{\ \ \ \ }\dfrac{\delta ^{2}p_{h}}{dt^{2}}+\{p_{s}R_{hji}^{\text{
\ \ \ \ }s}+\frac{1}{2}(\nabla _{i}c_{jh}+\nabla _{j}c_{ih}-\nabla
_{h}c_{ij})\}\dfrac{dx^{i}}{dt}\dfrac{dx^{j}}{dt}=0%
\end{array}%
\right.  \label{G.2}
\end{equation}%
where $\dfrac{\delta ^{2}p_{h}}{dt^{2}}=$ $\mathit{\ }\dfrac{d}{dt}(\mathit{%
\ }\dfrac{\delta p_{h}}{dt})-\Gamma _{ih}^{a}\mathit{\ \ }\dfrac{\delta p_{a}%
}{dt}\dfrac{dx^{i}}{dt}.$ Thus we have the following result.

\begin{theorem}
Let $\tilde{C}$ be a curve on $T^{\ast }M$ and locally expressed by $%
x^{h}=x^{h}(t),$ $x^{\bar{h}}=p_{h}(t)$ with respect to the induced
coordinates $(x^{h},x^{\overline{h}})$ in $\pi ^{-1}\left( U\right) \subset
T^{\ast }M$. The curve $\tilde{C}$ is a geodesic of the modified Riemannian
extension $\widetilde{g}_{\nabla ,c}$, if the projection $C$ of $\widetilde{C%
}$ is a geodesic in $M$ with the symmetric connection $\nabla $ and $p_{h}(t)
$ satisfies the differential equation (b) in (\ref{G.2}).\bigskip 
\end{theorem}

The Riemannian curvature tensor $\widetilde{R}$ of $T^{\ast }M$ with the
modified Riemannian extension $\widetilde{g}_{\nabla ,c}$ is obtained from
the well-known formula%
\begin{equation*}
\widetilde{R}\left( \widetilde{X},\widetilde{Y}\right) \widetilde{Z}=%
\widetilde{\nabla }_{\widetilde{X}}\widetilde{\nabla }_{\widetilde{Y}}%
\widetilde{Z}-\widetilde{\nabla }_{\widetilde{Y}}\widetilde{\nabla }_{%
\widetilde{X}}\widetilde{Z}-\widetilde{\nabla }_{\left[ \widetilde{X},%
\widetilde{Y}\right] }\widetilde{Z}
\end{equation*}%
for all $\widetilde{X},\widetilde{Y},\widetilde{Z}\in \Im _{0}^{1}(T^{\ast
}M)$. Then from Lemma \ref{Lemma1} and Proposition \ref{propo1}, we get the
following proposition.

\begin{proposition}
\label{propo2}The components of the curvature tensor $\widetilde{R}$ of the
cotangent bundle $T^{\ast }M$ with the modified Riemannian extension $%
\widetilde{g}_{\nabla ,c}$ are given as follows:%
\begin{eqnarray}
\widetilde{R}(E_{i},E_{j})E_{k} &=&R_{ijk}^{\text{ \ \ \ }h}E_{h}{\,}
\label{AB4.3} \\
&&{+\{p}_{s}(\nabla _{i}R_{hkj}^{\text{ \ \ \ }s}-\nabla _{j}R_{hki}^{\text{
\ \ \ }s})  \notag \\
&&+\frac{1}{2}\{\nabla _{i}(\nabla _{k}c_{jh}-\nabla _{h}c_{jk})-\nabla
_{j}(\nabla _{k}c_{ih}-\nabla _{h}c_{ik})  \notag \\
&&-R_{ijk}^{\text{ \ \ \ }m}c_{mh}-R_{ijh}^{\text{ \ \ \ }m}c_{km}\}\}E_{%
\overline{h}}  \notag \\
\widetilde{R}(E_{i},E_{j})E_{\overline{k}} &=&R_{jih}^{\text{ \ \ \ }k}E_{%
\overline{h}}  \notag \\
{\,}\widetilde{R}(E_{i},E_{\overline{j}})E_{k} &=&{\,-}R_{hki}^{\text{ \ \ \
\ }j}E_{\overline{h}}  \notag \\
{\,}\widetilde{R}(E_{\overline{i}},E_{j})E_{k} &=&R_{hkj}^{\text{ \ \ \ \ }%
i}E_{\overline{h}}  \notag \\
\widetilde{R}(E_{\overline{i}},E_{j})E_{\overline{k}} &=&0,{\,}\text{\ }%
\widetilde{R}(E_{i},E_{\overline{j}})E_{\overline{k}}={\,0,}  \notag \\
\widetilde{R}(E_{\overline{i}},E_{\overline{j}})E_{k} &=&0,\text{ }%
\widetilde{R}(E_{\overline{i}},E_{\overline{j}})E_{\overline{k}}=0  \notag
\end{eqnarray}%
with respect to the adapted frame $\left\{ E_{\alpha }\right\} .$
\end{proposition}

Proposition \ref{propo2} leads to the following result.

\begin{theorem}
Let $\nabla $ be a symmetric connection on $M$ and $T^{\ast }M$ be the
cotangent bundle with the modified Riemannian extension $\widetilde{g}%
_{\nabla ,c}$ over $(M,\nabla )$. Then $(T^{\ast }M,\widetilde{g}_{\nabla
,c})$ is locally flat if and only if $(M,\nabla )$ is locally flat and the
components $c_{ij}$ of $c$ satisfy the condition%
\begin{equation}
\nabla _{i}(\nabla _{k}c_{jh}-\nabla _{h}c_{jk})-\nabla _{j}(\nabla
_{k}c_{ih}-\nabla _{h}c_{ik})=0.  \label{AB4.4}
\end{equation}
\end{theorem}

Let $\widetilde{X}$ and $\widetilde{Y}$ be vector fields of $T^{\ast }M.$
The curvature operator $\widetilde{R}(\widetilde{X},\widetilde{Y})$ is a
differential operator on $T^{\ast }M.$ Similarly, for vector fields $X$ and $%
Y$ of $M$, $R(X,Y)$ is a differential operator on $M.$ Now, we operate the
curvature operator to the curvature tensor. That is, for all $\widetilde{Z},%
\widetilde{W}$ and $\widetilde{U}$, we consider the condition $(\widetilde{R}%
(\widetilde{X},\widetilde{Y})\widetilde{R})(\widetilde{Z},\widetilde{W})%
\widetilde{U}=0.$

The tensor $(\widetilde{R}(\widetilde{X},\widetilde{Y})\widetilde{R})(%
\widetilde{Z},\widetilde{W})\widetilde{U}$ has components%
\begin{eqnarray}
&&((\widetilde{R}(\widetilde{X},\widetilde{Y})\widetilde{R})(\widetilde{Z},%
\widetilde{W})\widetilde{U})_{\alpha \beta \gamma \theta \sigma }^{\text{ \
\ \ \ \ \ \ \ }\varepsilon }  \label{R.1} \\
&=&\widetilde{R}_{\alpha \beta \tau }^{\text{ \ \ \ }\varepsilon }\widetilde{%
R}_{\gamma \theta \sigma }^{\text{ \ \ \ }\tau }-\widetilde{R}_{\alpha \beta
\gamma }^{\text{ \ \ \ }\tau }\widetilde{R}_{\tau \theta \sigma }^{\text{ \
\ \ }\varepsilon }-\widetilde{R}_{\alpha \beta \theta }^{\text{ \ \ \ }\tau }%
\widetilde{R}_{\gamma \tau \sigma }^{\text{ \ \ \ }\varepsilon }-\widetilde{R%
}_{\alpha \beta \sigma }^{\text{ \ \ \ }\tau }\widetilde{R}_{\gamma \theta
\tau }^{\text{ \ \ \ }\varepsilon },  \notag
\end{eqnarray}%
with respect to the adapted frame $\left\{ E_{\alpha }\right\} $. Similarly,
for all $X,Y,Z,W,U$ on $M,$%
\begin{eqnarray*}
&&((R(X,Y)R)(Z,W)U)_{ijklm}^{\text{ \ \ \ \ \ \ \ }n} \\
&=&R_{ijp}^{\text{ \ \ \ }n}R_{klm}^{\text{ \ \ \ \ }p}-R_{ijk}^{\text{ \ \
\ }p}R_{plm}^{\text{ \ \ \ \ }n}-R_{ijl}^{\text{ \ \ \ }p}R_{kpm}^{\text{ \
\ \ \ }n}-R_{ijm}^{\text{ \ \ \ }p}R_{klp}^{\text{ \ \ \ \ }n} \\
&=&2\nabla _{\lbrack i}\nabla _{j]}R_{klm}^{\text{ \ \ \ \ }n}
\end{eqnarray*}%
where $2\nabla _{\lbrack i}\nabla _{j]}=\nabla _{i}\nabla _{j}-\nabla
_{j}\nabla _{i}.$

Case of $\alpha =i,\beta =j,\gamma =k,\theta =l,\sigma =\overline{m}%
,\varepsilon =\overline{n}$ in (\ref{R.1}), by virtue of (\ref{AB4.3}) the
equation (\ref{R.1}) reduces to%
\begin{eqnarray}
&&((\widetilde{R}(\widetilde{X},\widetilde{Y})\widetilde{R})(\widetilde{Z},%
\widetilde{W})\widetilde{U})_{ijkl\overline{m}}^{\text{ \ \ \ \ \ \ \ \ }%
\overline{n}}  \label{R.2} \\
&=&\widetilde{R}_{ijp}^{\text{ \ \ \ }\overline{n}}\widetilde{R}_{kl%
\overline{m}}^{\text{ \ \ \ \ }p}+\widetilde{R}_{ij\overline{p}}^{\text{ \ \
\ }\overline{n}}\widetilde{R}_{kl\overline{m}}^{\text{ \ \ \ \ }\overline{p}%
}-\widetilde{R}_{ijk}^{\text{ \ \ \ }p}\widetilde{R}_{pl\overline{m}}^{\text{
\ \ \ \ }\overline{n}}-\widetilde{R}_{ijk}^{\text{ \ \ \ }\overline{p}}%
\widetilde{R}_{\overline{p}l\overline{m}}^{\text{ \ \ \ \ }\overline{n}} 
\notag \\
&&-\widetilde{R}_{ijl}^{\text{ \ \ \ }p}\widetilde{R}_{kp\overline{m}}^{%
\text{ \ \ \ \ }\overline{n}}-\widetilde{R}_{ijl}^{\text{ \ \ \ }\overline{p}%
}\widetilde{R}_{k\overline{p}\overline{m}}^{\text{ \ \ \ \ }\overline{n}}-%
\widetilde{R}_{ij\overline{m}}^{\text{ \ \ \ }p}\widetilde{R}_{klp}^{\text{
\ \ \ \ }\overline{n}}-\widetilde{R}_{ij\overline{m}}^{\text{ \ \ \ }%
\overline{p}}\widetilde{R}_{kl\overline{p}}^{\text{ \ \ \ \ }\overline{n}} 
\notag \\
&=&-(R_{ijp}^{\text{ \ \ \ }m}R_{kl\text{ }n}^{\text{ \ \ \ \ }p}-R_{ijk}^{%
\text{ \ \ \ }p}R_{pl\text{ }n}^{\text{ \ \ \ \ }m}-R_{ijl}^{\text{ \ \ \ }%
p}R_{kpn}^{\text{ \ \ \ \ }m}-R_{ijn}^{\text{ \ \ \ }p}R_{klp}^{\text{ \ \ \
\ }m}  \notag \\
&=&-2\nabla _{\lbrack i}\nabla _{j]}R_{kl\text{ }n}^{\text{ \ \ \ \ \ }m} 
\notag \\
&=&-((R(X,Y)R)(Z,W)U)_{ijkl\text{ }n}^{\text{ \ \ \ \ \ \ \ }m}.  \notag
\end{eqnarray}

Case of $\alpha =i,\beta =j,\gamma =\overline{k},\theta =l,\sigma
=m,\varepsilon =\overline{n}$ in (\ref{R.1}), by virtue of (\ref{AB4.3}), we
have%
\begin{eqnarray}
&&((\widetilde{R}(\widetilde{X},\widetilde{Y})\widetilde{R})(\widetilde{Z},%
\widetilde{W})\widetilde{U})_{ij\overline{k}lm}^{\text{ \ \ \ \ \ \ \ \ }%
\overline{n}}  \label{R.3} \\
&=&\widetilde{R}_{ijp}^{\text{ \ \ \ }\overline{n}}\widetilde{R}_{\overline{k%
}lm}^{\text{ \ \ \ \ }p}+\widetilde{R}_{ij\overline{p}}^{\text{ \ \ \ }%
\overline{n}}\widetilde{R}_{\overline{k}lm}^{\text{ \ \ \ \ }\overline{p}}-%
\widetilde{R}_{ij\overline{k}}^{\text{ \ \ \ }p}\widetilde{R}_{plm}^{\text{
\ \ \ \ }\overline{n}}-\widetilde{R}_{ij\overline{k}}^{\text{ \ \ \ }%
\overline{p}}\widetilde{R}_{\overline{p}lm}^{\text{ \ \ \ \ }\overline{n}} 
\notag \\
&&-\widetilde{R}_{ijl}^{\text{ \ \ \ }p}\widetilde{R}_{\overline{k}pm}^{%
\text{ \ \ \ \ }\overline{n}}-\widetilde{R}_{ijl}^{\text{ \ \ \ }\overline{p}%
}\widetilde{R}_{\overline{k}\overline{p}m}^{\text{ \ \ \ \ }\overline{n}}-%
\widetilde{R}_{ijm}^{\text{ \ \ \ }p}\widetilde{R}_{\overline{k}lp}^{\text{
\ \ \ \ }\overline{n}}-\widetilde{R}_{ijm}^{\text{ \ \ \ }\overline{p}}%
\widetilde{R}_{\overline{k}l\overline{p}}^{\text{ \ \ \ \ }\overline{n}} 
\notag \\
&=&-(R_{ijp}^{\text{ \ \ \ }k}R_{nlm}^{\text{ \ \ \ \ }p}-R_{ijn}^{\text{ \
\ \ }p}R_{plm}^{\text{ \ \ \ \ }k}-R_{ijl}^{\text{ \ \ \ }p}R_{npm}^{\text{
\ \ \ \ }k}-R_{ijm}^{\text{ \ \ \ }p}R_{nlp}^{\text{ \ \ \ \ }k}  \notag \\
&=&-2\nabla _{\lbrack i}\nabla _{j]}R_{nlm}^{\text{ \ \ \ \ }k}  \notag \\
&=&-((R(X,Y)R)(Z,W)U)_{ijnlm}^{\text{ \ \ \ \ \ \ \ }k}.  \notag
\end{eqnarray}

Case of $\alpha =i,\beta =j,\gamma =\overline{k},\theta =\overline{l},\sigma
=\overline{m},\varepsilon =\overline{n}$ in (\ref{R.1}), by virtue of (\ref%
{AB4.3}), it follows that 
\begin{equation}
((\widetilde{R}(\widetilde{X},\widetilde{Y})\widetilde{R})(\widetilde{Z},%
\widetilde{W})\widetilde{U})_{ij\overline{k}\overline{l}\overline{m}}^{\text{
\ \ \ \ \ \ \ \ }\overline{n}}=0.  \label{R.4}
\end{equation}

Case of $\alpha =i,\beta =j,\gamma =k,\theta =l,\sigma =m,\varepsilon =%
\overline{n}$ in (\ref{R.1}), we obtain%
\begin{eqnarray}
&&((\widetilde{R}(\widetilde{X},\widetilde{Y})\widetilde{R})(\widetilde{Z},%
\widetilde{W})\widetilde{U})_{ijklm}^{\text{ \ \ \ \ \ \ \ \ }\overline{n}}
\label{R.5} \\
&=&\widetilde{R}_{ijp}^{\text{ \ \ \ }\overline{n}}\widetilde{R}_{klm}^{%
\text{ \ \ \ \ }p}+\widetilde{R}_{ij\overline{p}}^{\text{ \ \ \ }\overline{n}%
}\widetilde{R}_{klm}^{\text{ \ \ \ \ }\overline{p}}-\widetilde{R}_{ijk}^{%
\text{ \ \ \ }p}\widetilde{R}_{plm}^{\text{ \ \ \ \ }\overline{n}}-%
\widetilde{R}_{ijk}^{\text{ \ \ \ }\overline{p}}\widetilde{R}_{\overline{p}%
lm}^{\text{ \ \ \ \ }\overline{n}}  \notag \\
&&-\widetilde{R}_{ijl}^{\text{ \ \ \ }p}\widetilde{R}_{kpm}^{\text{ \ \ \ \ }%
\overline{n}}-\widetilde{R}_{ijl}^{\text{ \ \ \ }\overline{p}}\widetilde{R}%
_{k\overline{p}m}^{\text{ \ \ \ \ }\overline{n}}-\widetilde{R}_{ijm}^{\text{
\ \ \ }p}\widetilde{R}_{klp}^{\text{ \ \ \ \ }\overline{n}}-\widetilde{R}%
_{ijm}^{\text{ \ \ \ }\overline{p}}\widetilde{R}_{kl\overline{p}}^{\text{ \
\ \ \ }\overline{n}}.  \notag
\end{eqnarray}

Case of $\alpha =\overline{i},\beta =j,\gamma =k,\theta =l,\sigma
=m,\varepsilon =\overline{n}$ in (\ref{R.1}), by virtue of (\ref{AB4.3}), we
get%
\begin{eqnarray}
&&((\widetilde{R}(\widetilde{X},\widetilde{Y})\widetilde{R})(\widetilde{Z},%
\widetilde{W})\widetilde{U})_{\overline{i}jklm}^{\text{ \ \ \ \ \ \ \ \ }%
\overline{n}}  \label{R.6} \\
&=&\widetilde{R}_{\overline{i}jp}^{\text{ \ \ \ }\overline{n}}\widetilde{R}%
_{klm}^{\text{ \ \ \ \ }p}+\widetilde{R}_{\overline{i}j\overline{p}}^{\text{
\ \ \ }\overline{n}}\widetilde{R}_{klm}^{\text{ \ \ \ \ }\overline{p}}-%
\widetilde{R}_{\overline{i}jk}^{\text{ \ \ \ }p}\widetilde{R}_{plm}^{\text{
\ \ \ \ }\overline{n}}-\widetilde{R}_{\overline{i}jk}^{\text{ \ \ \ }%
\overline{p}}\widetilde{R}_{\overline{p}lm}^{\text{ \ \ \ \ }\overline{n}} 
\notag \\
&&-\widetilde{R}_{\overline{i}jl}^{\text{ \ \ \ }p}\widetilde{R}_{kpm}^{%
\text{ \ \ \ \ }\overline{n}}-\widetilde{R}_{\overline{i}jl}^{\text{ \ \ \ }%
\overline{p}}\widetilde{R}_{k\overline{p}m}^{\text{ \ \ \ \ }\overline{n}}-%
\widetilde{R}_{\overline{i}jm}^{\text{ \ \ \ }p}\widetilde{R}_{klp}^{\text{
\ \ \ \ }\overline{n}}-\widetilde{R}_{\overline{i}jm}^{\text{ \ \ \ }%
\overline{p}}\widetilde{R}_{kl\overline{p}}^{\text{ \ \ \ \ }\overline{n}} 
\notag \\
&=&R_{npj}^{\text{ \ \ \ }i}R_{klm}^{\text{ \ \ \ \ }p}-R_{pkj}^{\text{ \ \
\ \ }i}R_{nml}^{\text{ \ \ \ \ }p}+R_{plj}^{\text{ \ \ \ \ }i}R_{nmk}^{\text{
\ \ \ \ }p}-R_{pmj}^{\text{ \ \ \ \ }i}R_{lkn}^{\text{ \ \ \ \ }p}  \notag \\
&=&2\nabla _{\lbrack n}\nabla _{m]}R_{lkj}^{\text{ \ \ \ \ }i}-2\nabla
_{\lbrack k}\nabla _{l]}R_{nmj}^{\text{ \ \ \ \ }i}  \notag \\
&=&((R(X,Y)R)(Z,W)U)_{nmlkj}^{\text{ \ \ \ \ \ \ \ }i}-((R(X,Y)R)(Z,W)U)_{kl%
\text{ }nmj}^{\text{ \ \ \ \ \ \ \ }k}.  \notag
\end{eqnarray}

The other coefficients of $(\widetilde{R}(\widetilde{X},\widetilde{Y})%
\widetilde{R})(\widetilde{Z},\widetilde{W})\widetilde{U}$ reduce to one of (%
\ref{R.2}), (\ref{R.3}) or (\ref{R.4}) by the property of the curvature
tensor.

From (\ref{R.2})-(\ref{R.6}), we have the following result.

\begin{theorem}
\label{Theo1}Let $\nabla $ be a symmetric connection on $M$ and $T^{\ast }M$
be the cotangent bundle with the modified Riemannian extension $\widetilde{g}%
_{\nabla ,c}$ over $(M,\nabla )$. Then $(\widetilde{R}(\widetilde{X},%
\widetilde{Y})\widetilde{R})(\widetilde{Z},\widetilde{W})\widetilde{U}=0$ if
and only if the following conditions hold:

i)$(R(X,Y)R)(Z,W)U=0$

ii) $R_{ijk}^{\text{ \ \ \ }\overline{h}}=0,$ from which it follows that $%
\nabla _{i}R_{hkj}^{\text{ \ \ \ }s}-\nabla _{j}R_{hki}^{\text{ \ \ \ }s}=0$
and $\nabla _{i}(\nabla _{k}c_{jh}-\nabla _{h}c_{jk})-\nabla _{j}(\nabla
_{k}c_{ih}-\nabla _{h}c_{ik})-R_{ijk}^{\text{ \ \ \ }m}c_{mh}-R_{ijh}^{\text{
\ \ \ }m}c_{km}=0,$ where $R$ and $\widetilde{R}$ is the curvature tensors
of the symmetric connection $\nabla $ and the Levi-Civita connection $%
\widetilde{\nabla }$ of the modified Riemannian extension $\widetilde{g}%
_{\nabla ,c}$ respectively.
\end{theorem}

Theorem \ref{Theo1} immediately give the following result.

\begin{corollary}
If the symmetric connection $\nabla $ on $M$ is locally symmetric, then $(%
\widetilde{R}(\widetilde{X},\widetilde{Y})\widetilde{R})(\widetilde{Z},%
\widetilde{W})\widetilde{U}=0$ if and only if 
\begin{equation}
\nabla _{i}(\nabla _{k}c_{jh}-\nabla _{h}c_{jk})-\nabla _{j}(\nabla
_{k}c_{ih}-\nabla _{h}c_{ik})-R_{ijk}^{\text{ \ \ \ }m}c_{mh}-R_{ijh}^{\text{
\ \ \ }m}c_{km}=0.  \label{R.7}
\end{equation}
\end{corollary}

Now, we consider the components of $\widetilde{\nabla }\widetilde{R}$. Using
Proposition \ref{propo1} and (\ref{AB4.3}), by direct computation, we obtain
following relations%
\begin{eqnarray*}
\widetilde{\nabla }_{l}\widetilde{R}_{ijk}^{\text{ \ \ \ }h} &=&\nabla
_{l}R_{ijk}^{\text{ \ \ \ }h}, \\
\widetilde{\nabla }_{l}\widetilde{R}_{ijk}^{\text{ \ \ \ }\overline{h}} &=&{p%
}_{s}(\nabla _{l}\nabla _{i}R_{hkj}^{\text{ \ \ \ }s}-\nabla _{l}\nabla
_{j}R_{hki}^{\text{ \ \ \ }s})+\frac{1}{2}\{\nabla _{l}\nabla _{i}(\nabla
_{k}c_{jh}-\nabla _{h}c_{jk}) \\
&&-\nabla _{l}\nabla _{j}(\nabla _{k}c_{ih}-\nabla _{h}c_{ik})-(\nabla
_{l}R_{ijk}^{\text{ \ \ \ }m})c_{mh}-R_{ijk}^{\text{ \ \ \ }m}(\nabla
_{l}c_{mh}) \\
&&-(\nabla _{l}R_{ijh}^{\text{ \ \ \ }m})c_{km}-R_{ijh}^{\text{ \ \ \ }%
m}(\nabla _{l}c_{km})\} \\
\widetilde{\nabla }_{l}\widetilde{R}_{ij\overline{k}}^{\text{ \ \ \ }%
\overline{h}} &=&\nabla _{l}R_{jih}^{\text{ \ \ \ }k}, \\
\widetilde{\nabla }_{l}\widetilde{R}_{i\overline{j}k}^{\text{ \ \ \ }%
\overline{h}} &=&-\nabla _{l}R_{hki}^{\text{ \ \ \ }j}, \\
\widetilde{\nabla }_{l}\widetilde{R}_{\overline{i}jk}^{\text{ \ \ \ }%
\overline{h}} &=&\nabla _{l}R_{hkj}^{\text{ \ \ \ \ }i}, \\
\widetilde{\nabla }_{\overline{l}}\widetilde{R}_{ijk}^{\text{ \ \ \ }%
\overline{h}} &=&\nabla _{i}R_{hkj}^{\text{ \ \ \ }l}-\nabla _{j}R_{hki}^{%
\text{ \ \ \ }l},
\end{eqnarray*}%
all the others being zero, with respect to the adapted frame $\left\{
E_{\alpha }\right\} .$ Hence we have the following.

\begin{theorem}
Let $\nabla $ be a symmetric connection on $M$ and $T^{\ast }M$ be the
cotangent bundle with the modified Riemannian extension $\widetilde{g}%
_{\nabla ,c}$ over $(M,\nabla )$. Then $(T^{\ast }M,\widetilde{g}_{\nabla
,c})$ is locally symmetric if and only if $(M,\nabla )$ is locally symmetric
and the components $c_{ij}$ of $c$ satisfy the condition%
\begin{eqnarray}
&&\nabla _{l}\nabla _{i}(\nabla _{k}c_{jh}-\nabla _{h}c_{jk})-\nabla
_{l}\nabla _{j}(\nabla _{k}c_{ih}-\nabla _{h}c_{ik})  \label{AB4.5} \\
&&-R_{ijk}^{\text{ \ \ \ }m}(\nabla _{l}c_{mh})-R_{ijh}^{\text{ \ \ \ }%
m}(\nabla _{l}c_{km})  \notag \\
&=&0.  \notag
\end{eqnarray}
\end{theorem}

We now turn our attention to the Ricci tensor and scalar curvature of the
modified Riemannian extension $\widetilde{g}_{\nabla ,c}$. Let $\widetilde{R}%
_{\alpha \beta }=$ $\widetilde{R}_{\sigma \alpha \beta }^{\text{ \ \ \ \ \ \ 
}\sigma }$ and $\widetilde{r}=(\widetilde{g}_{\nabla ,c})^{\alpha \beta }$ $%
\widetilde{R}_{\alpha \beta }$ denote the Ricci tensor and scalar curvature
of the modified Riemannian extension $\widetilde{g}_{\nabla ,c}$,
respectively. From (\ref{AB4.3}), the components of the Ricci tensor $%
R_{\alpha \beta }$ are characterized by%
\begin{eqnarray}
{}\widetilde{{R}}{_{jk}}{=R_{jk}+} &&R_{kj}  \label{RC} \\
\widetilde{{R}}{_{\overline{j}k}}{=0,} &&  \notag \\
\text{ }\widetilde{{R}}{_{j\overline{k}}}{=0,} &&  \notag \\
\text{ }\widetilde{{R}}{_{\overline{j}\overline{k}}}{=0,} &&  \notag
\end{eqnarray}%
with respect to the adapted frame $\left\{ E\alpha \right\} $. \ Next we have

\begin{theorem}
Let $\nabla $ be a symmetric connection on $M$ and $T^{\ast }M$ be the
cotangent bundle with the modified Riemannian extension $\widetilde{g}%
_{\nabla ,c}$ over $(M,\nabla )$. Then $(T^{\ast }M,\widetilde{g}_{\nabla
,c})$ is Ricci flat if and only if the Ricci tensor of $\nabla $ is skew
symmetric (for Riemannian extension, see \cite{Patterson}).
\end{theorem}

Now, we operate the curvature operator to the Ricci tensor. The tensors $(%
\widetilde{R}(\widetilde{X},\widetilde{Y})\widetilde{Ric})(\widetilde{Z},%
\widetilde{W})$ and $(R(X,Y)Ric)(Z,W)$ have coefficients%
\begin{equation*}
((\widetilde{R}(\widetilde{X},\widetilde{Y})\widetilde{Ric})(\widetilde{Z},%
\widetilde{W}))_{\alpha \beta \gamma \theta }=\widetilde{R}_{\alpha \beta
\gamma }^{\text{ \ \ \ }\varepsilon }\widetilde{R}_{\varepsilon \theta }+%
\widetilde{R}_{\alpha \beta \theta }^{\text{ \ \ \ }\varepsilon }\widetilde{R%
}_{\gamma \varepsilon }
\end{equation*}%
and%
\begin{equation*}
((R(X,Y)Ric)(Z,W))_{ijkl}=R_{ijk}^{\text{ \ \ \ }p}R_{pl}+R_{ijl}^{\text{ \
\ \ }p}R_{kp}
\end{equation*}%
respectively. By putting $\alpha =i,\beta =j,\gamma =k,\theta =l$, it
follows that%
\begin{eqnarray*}
((\widetilde{R}(\widetilde{X},\widetilde{Y})\widetilde{Ric})(\widetilde{Z},%
\widetilde{W}))_{ijkl} &=&\widetilde{R}_{ijk}^{\text{ \ \ \ }p}\widetilde{R}%
_{pl}+\widetilde{R}_{ijl}^{\text{ \ \ \ }p}\widetilde{R}_{kp} \\
&=&R_{ijk}^{\text{ \ \ \ }p}(R_{pl}+R_{lp})+R_{ijl}^{\text{ \ \ \ }%
p}(R_{kp}+R_{pk})
\end{eqnarray*}%
all the others being zero. Let the base manifold $M$ be a Riemannian
manifold with the metric $g$ and $\nabla $ be the Levi-Civita connection of $%
g$. Then 
\begin{eqnarray*}
((\widetilde{R}(\widetilde{X},\widetilde{Y})\widetilde{Ric})(\widetilde{Z},%
\widetilde{W}))_{ijkl} &=&2R_{ijk}^{\text{ \ \ \ }p}R_{pl}+2R_{ijl}^{\text{
\ \ \ }p}R_{kp} \\
&=&2((R(X,Y)Ric)(Z,W))_{ijkl}.
\end{eqnarray*}%
Therefore we get the following.

\begin{theorem}
Let $\nabla $ be the Levi-Civita connection on the Riemannian manifold $(M,g)
$ and $T^{\ast }M$ be its cotangent bundle with the modified Riemannian
extension $\widetilde{g}_{\nabla ,c}$. Then $(\widetilde{R}(\widetilde{X},%
\widetilde{Y})\widetilde{Ric})(\widetilde{Z},\widetilde{W})=0$ if and only
if $(R(X,Y)Ric)(Z,W)=0.$
\end{theorem}

From (\ref{AB4.2}) and (\ref{RC}), the scalar curvature of the modified
Riemannian extension $\widetilde{g}_{\nabla ,c}$ is given by%
\begin{equation*}
\widetilde{r}=0.
\end{equation*}

\begin{theorem}
Let $\nabla $ be a symmetric connection on $M$ and $T^{\ast }M$ be the
cotangent bundle with the modified Riemannian extension $\widetilde{g}%
_{\nabla ,c}$ over $(M,\nabla )$. Then $(T^{\ast }M,\widetilde{g}_{\nabla
,c})$ is a space of constant scalar curvature $0.$
\end{theorem}

In the following we give the conditions under which the cotangent bundle $%
(T^{\ast }M,\widetilde{g}_{\nabla ,c})$ is locally conformally flat. The
cotangent bundle $T^{\ast }M$ with the modified Riemannian extension $%
\widetilde{g}_{\nabla ,c}$ is locally conformally flat if and only if its
Weyl tensor $\widetilde{W}$ vanishes, where the Weyl tensor is given by%
\begin{eqnarray*}
\widetilde{W}_{\alpha \beta \gamma \sigma } &=&\widetilde{{R}}_{\alpha \beta
\gamma \sigma }+{\frac{\widetilde{r}}{2(2n-1)(n-1)}\left\{ (\widetilde{g}%
_{\nabla ,c})_{\alpha \gamma }(\widetilde{g}_{\nabla ,c})_{\beta \sigma }-(%
\widetilde{g}_{\nabla ,c})_{\alpha \sigma }(\widetilde{g}_{\nabla
,c})_{\beta \gamma }\right\} } \\
&&-{\frac{1}{2(n-1)}(}(\widetilde{g}_{\nabla ,c})_{\beta \sigma }\widetilde{{%
R}}_{\alpha \gamma }{-}(\widetilde{g}_{\nabla ,c})_{\alpha \sigma }%
\widetilde{{R}}_{\beta \gamma }{+}(\widetilde{g}_{\nabla ,c})_{\alpha \gamma
}\widetilde{{R}}_{\beta \sigma }{-}(\widetilde{g}_{\nabla ,c})_{\beta \gamma
}\widetilde{{R}}_{\alpha \sigma }{),}
\end{eqnarray*}%
where $\widetilde{{R}}_{\alpha \beta \gamma \sigma }=\widetilde{{R}}_{\alpha
\beta \gamma }^{\text{ \ \ \ \ }\lambda }(\widetilde{g}_{\nabla
,c})_{\lambda \sigma }.$

From \ref{AB4.3}, we obtain 
\begin{eqnarray*}
\widetilde{R}_{ijkn} &=&R_{ijk}^{\text{ \ \ }h}c_{hn}+p_{s}(\nabla
_{i}R_{nkj}^{\text{ \ \ \ }s}-\nabla _{j}R_{nki}^{\text{ \ \ \ }s}) \\
&&+\frac{1}{2}\{\nabla _{i}(\nabla _{k}c_{jn}-\nabla _{n}c_{jk})-\nabla
_{j}(\nabla _{k}c_{in}-\nabla _{n}c_{ik})-R_{ijk}^{\text{ \ \ }%
h}c_{hn}-R_{ijn}^{\text{ \ \ }h}c_{kh}\}
\end{eqnarray*}%
\begin{equation*}
\widetilde{R}_{ijk\overline{n}}=R_{ijk}^{\text{ \ \ \ }n}
\end{equation*}%
\begin{equation*}
\widetilde{R}_{ij\overline{k}n}=R_{jin}^{\text{ \ \ \ }k}
\end{equation*}%
\begin{equation*}
\widetilde{R}_{i\overline{j}kn}=R_{kni}^{\text{ \ \ \ }j}
\end{equation*}%
\begin{equation*}
\widetilde{R}_{\overline{i}jkn}=R_{nkj}^{\text{ \ \ \ }i}
\end{equation*}%
\begin{equation*}
otherwise=0.
\end{equation*}

The non-zero components of Weyl tensor of the modified Riemannian extension $%
\widetilde{g}_{\nabla ,c}$ are given by%
\begin{eqnarray*}
\widetilde{W}_{ijkn} &=&\widetilde{R}_{ijkn}-\frac{1}{2(n-1)}%
\{c_{jn}(R_{ik}+R_{ki})-c_{in}(R_{jk}+R_{kj}) \\
&&+c_{ik}(R_{jn}+R_{nj})-c_{jk}(R_{in}+R_{ni})\}
\end{eqnarray*}%
\begin{eqnarray*}
\widetilde{W}_{ijk\overline{n}} &=&\widetilde{R}_{ijk\overline{n}}-\frac{1}{%
2(n-1)}(\delta _{j}^{n}(R_{ik}+R_{ki})-\delta _{i}^{n}(R_{jk}+R_{kj})) \\
\widetilde{W}_{ij\overline{k}n} &=&\widetilde{R}_{ij\overline{k}n}-\frac{1}{%
2(n-1)}(\delta _{i}^{k}(R_{jn}+R_{nj})-\delta _{j}^{k}(R_{in}+R_{ni})) \\
\widetilde{W}_{i\overline{j}kn} &=&\widetilde{R}_{i\overline{j}kn}-\frac{1}{%
2(n-1)}(\delta _{n}^{j}(R_{ik}+R_{ki})-\delta _{j}^{k}(R_{in}+R_{ni})) \\
\widetilde{W}_{\overline{i}jkn} &=&\widetilde{R}_{\overline{i}jkn}-\frac{1}{%
2(n-1)}(\delta _{k}^{i}(R_{jn}+R_{nj})-\delta _{n}^{i}(R_{jk}+R_{kj}))
\end{eqnarray*}%
Using the same way in \cite{Afifi}, we can say the same result in \cite%
{Afifi} by means of adapted frame.

\begin{theorem}
\bigskip Let $\nabla $ be a symmetric connection on $M$ and $T^{\ast }M$ be
the cotangent bundle with the modified Riemannian extension $\widetilde{g}%
_{\nabla ,c}$ over $(M,\nabla )$. Then $(T^{\ast }M,\widetilde{g}_{\nabla
,c})$ is locally conformally flat if and only if $(M,\nabla )$ is
projectively flat and the components $c_{ij}$ of $c$ satisfy the condition%
\begin{equation}
\nabla _{i}(\nabla _{k}c_{jn}-\nabla _{n}c_{jk})-\nabla _{j}(\nabla
_{k}c_{in}-\nabla _{n}c_{ik})-R_{ijk}^{\text{ \ \ }h}c_{hn}-R_{ijn}^{\text{
\ \ }h}c_{kh}=0\text{ (also see, \cite{Afifi})}.  \label{AB4.6}
\end{equation}
\end{theorem}

Finally, we finish the section by the projective curvature tensor of the
modified Riemannian extension $\widetilde{g}_{\nabla ,c}.$ A manifold is
said to be projectively flat if the projective curvature tensor vanishes.
The projective curvature tensor is defined by%
\begin{equation*}
\widetilde{P}_{\alpha \beta \gamma \sigma }=\widetilde{{R}}_{\alpha \beta
\gamma \sigma }-{\frac{1}{(2n-1)}(}(\widetilde{g}_{\nabla ,c})_{\alpha
\sigma }\widetilde{{R}}_{\beta \gamma }-(\widetilde{g}_{\nabla ,c})_{\beta
\sigma }\widetilde{{R}}_{\alpha \gamma }).
\end{equation*}

Case of $\alpha =i,\beta =j,\gamma =\overline{k},\sigma =n,$ we have%
\begin{eqnarray*}
\widetilde{P}_{ij\overline{k}n} &=&\widetilde{{R}}_{ij\overline{k}n}-{\frac{1%
}{(2n-1)}(}(\widetilde{g}_{\nabla ,c})_{in}\widetilde{{R}}_{j\overline{k}}-(%
\widetilde{g}_{\nabla ,c})_{jn}\widetilde{{R}}_{i\overline{k}}) \\
&=&R_{jin}^{\text{ \ \ \ }k}.
\end{eqnarray*}

Case of $\alpha =i,\beta =j,\gamma =k,\sigma =n,$ we obtain%
\begin{eqnarray*}
\widetilde{P}_{ijkn} &=&\widetilde{{R}}_{ijkn}-{\frac{1}{(2n-1)}(}(%
\widetilde{g}_{\nabla ,c})_{in}\widetilde{{R}}_{jk}-(\widetilde{g}_{\nabla
,c})_{jn}\widetilde{{R}}_{ik}) \\
&=&R_{ijk}^{\text{ \ \ }h}c_{hn}+p_{s}(\nabla _{i}R_{nkj}^{\text{ \ \ \ }%
s}-\nabla _{j}R_{nki}^{\text{ \ \ \ }s}) \\
&&+\frac{1}{2}\{\nabla _{i}(\nabla _{k}c_{jn}-\nabla _{n}c_{jk})-\nabla
_{j}(\nabla _{k}c_{in}-\nabla _{n}c_{ik})-R_{ijk}^{\text{ \ \ }%
h}c_{hn}-R_{ijn}^{\text{ \ \ }h}c_{kh}\} \\
&&-{\frac{1}{(2n-1)}(}c_{in}{(R}_{jk}+R_{kj})-c_{jn}({R}_{ik}+R_{ki})).
\end{eqnarray*}

The above equations give the following result.

\begin{theorem}
\bigskip Let $\nabla $ be a symmetric connection on $M$ and $T^{\ast }M$ be
the cotangent bundle with the modified Riemannian extension $\widetilde{g}%
_{\nabla ,c}$ over $(M,\nabla )$. Then $(T^{\ast }M,\widetilde{g}_{\nabla
,c})$ is projectively flat if and only if $(M,\nabla )$ is flat and the
components $c_{ij}$ of $c$ satisfy the condition%
\begin{equation}
\nabla _{i}(\nabla _{k}c_{jn}-\nabla _{n}c_{jk})-\nabla _{j}(\nabla
_{k}c_{in}-\nabla _{n}c_{ik})=0.  \label{AB4.7}
\end{equation}

\begin{remark}
i) If $c_{ij}=0,$ then the conditions (\ref{AB4.4}), (\ref{R.7}), (\ref%
{AB4.5}), (\ref{AB4.6}) and (\ref{AB4.7}) are identically fulfilled.

ii) If $c_{ij}$ is parallel with respect to $\nabla $, then the conditions (%
\ref{AB4.4}), (\ref{R.7}), (\ref{AB4.5}),\ (\ref{AB4.6}) and (\ref{AB4.7})
are identically fulfilled.

iii) If $c_{ij}$ satisfies the relation $\nabla _{i}c_{jk}-\nabla
_{j}c_{ik}=\nabla _{k}\omega _{ij}$, where the components $\omega _{ij}$
define a $2$--form on $M$ and if $(M,\nabla )$ is flat then the condition (%
\ref{AB4.4}), (\ref{R.7}), (\ref{AB4.5}), (\ref{AB4.6}) and (\ref{AB4.7})
are identically verified.
\end{remark}
\end{theorem}

\section{Curvature properties of another metric connection of the modified
Riemannian extension $\protect\widetilde{g}_{\protect\nabla ,c}$}

Let $\nabla $ be a linear connection on an $n-$dimensional differentiable
manifold $M$. The connection $\nabla $ is symmetric if its torsion tensor
vanishes, otherwise it is non-symmetric. If there is a Riemannian metric $g$
on $M$ such that $\nabla g=0$, then the connection $\nabla $ is a metric
connection, otherwise it is non-metric. It is well known that a linear
connection is symmetric and metric if and only if it is the Levi-Civita
connection. In section 4, we have considered the Levi-Civita connection $%
\widetilde{\nabla }$ of the modified Riemannian extension $\widetilde{g}%
_{\nabla ,c}$ on the cotangent bundle $T^{\ast }M$ over $(M,\nabla ).$ The
connection is the unique connection which satisfies $\widetilde{\nabla }%
_{\alpha }(\widetilde{g}_{\nabla ,c})_{\beta \gamma }=0$ and has a zero
torsion. Hayden \cite{Hayden} introduced a metric connection with a non-zero
torsion on a Riemannian manifold. Now we are interested in a metric
connection $^{(M)}\widetilde{\nabla }$ of the modified Riemannian extension $%
\widetilde{g}_{\nabla ,c}$ whose torsion tensor $^{^{(M)}\widetilde{\nabla }%
}T_{\gamma \beta }^{\varepsilon }$ is skew-symmetric in the indices $\gamma $
and $\beta .$ We denote components of the connection $^{(M)}\widetilde{%
\nabla }$ by $^{(M)}\widetilde{\Gamma }.$ The metric connection $^{(M)}%
\widetilde{\nabla }$ satisfies 
\begin{equation}
^{(M)}\widetilde{\nabla }_{\alpha }(\widetilde{g}_{\nabla ,c})_{\beta \gamma
}=0\text{ and }^{(M)}\widetilde{\Gamma }_{\alpha \beta }^{\gamma }-\text{ }%
^{(M)}\widetilde{\Gamma }_{\beta \alpha }^{\gamma }=\text{ }^{^{(M)}%
\widetilde{\nabla }}T_{\alpha \beta }^{\gamma }.  \label{AB5.1}
\end{equation}%
On the equation (\ref{AB5.1}) is solved with respect to $^{(M)}\widetilde{%
\Gamma }_{\alpha \beta }^{\gamma },$ one finds the following solution \cite%
{Hayden}%
\begin{equation}
^{(M)}\widetilde{\Gamma }_{\alpha \beta }^{\gamma }=\widetilde{\Gamma }%
_{\alpha \beta }^{\gamma }+\widetilde{U}_{\alpha \beta }^{\gamma },
\label{AB5.2}
\end{equation}%
where $\widetilde{\Gamma }_{\alpha \beta }^{\gamma }$ is components of the
Levi-Civita connection of the modified Riemannian extension $\widetilde{g}%
_{\nabla ,c}$, 
\begin{equation}
\widetilde{U}_{\alpha \beta \gamma }=\frac{1}{2}(^{^{(M)}\widetilde{\nabla }%
}T_{\alpha \beta \gamma }+\text{ }^{^{(M)}\widetilde{\nabla }}T_{\gamma
\alpha \beta }+\text{ }^{^{(M)}\widetilde{\nabla }}T_{\gamma \beta \alpha })
\label{AB5.3}
\end{equation}%
and%
\begin{equation*}
\widetilde{U}_{\alpha \beta \gamma }=\widetilde{U}_{\alpha \beta }^{\epsilon
}(\widetilde{g}_{\nabla ,c})_{\epsilon \gamma },\text{ }^{^{(M)}\widetilde{%
\nabla }}T_{\alpha \beta \gamma }=T_{\alpha \beta }^{\epsilon }(\widetilde{g}%
_{\nabla ,c})_{\epsilon \gamma }.
\end{equation*}

If we put%
\begin{equation}
^{^{(M)}\widetilde{\nabla }}T_{ij}^{\overline{r}}=-p_{s}R_{ijr}^{\text{ \ \ }%
s}  \label{AB5.4}
\end{equation}%
all other $^{^{(M)}\widetilde{\nabla }}T_{\alpha \beta }^{\gamma }$ not
related to $^{^{(M)}\widetilde{\nabla }}T_{ij}^{\overline{r}}$ being assumed
to be zero. \ We choose this $^{^{(M)}\widetilde{\nabla }}T_{\alpha \beta
}^{\gamma }$ in $T^{\ast }M$ which is skew-symmetric in the indices $\gamma $
and $\beta $ as torsion tensor and determine a metric connection in $T^{\ast
}M$ with respect to the modified Riemannian extension $\widetilde{g}_{\nabla
,c}.$ By using (\ref{AB5.3}) and (\ref{AB5.4}), we get non-zero component of 
$\widetilde{U}_{\alpha \beta }^{\gamma }$ as follows: 
\begin{equation*}
\widetilde{{U}}{_{ij}^{\overline{h}}}{=}p_{s}R_{jhi}^{\text{ \ \ }s}
\end{equation*}%
with respect to the adapted frame. From (\ref{AB5.2}) and Proposition \ref%
{propo1}, we have

\begin{proposition}
Let $\nabla $ be a symmetric connection on $M$ and $T^{\ast }M$ be the
cotangent bundle with the modified Riemannian extension $\widetilde{g}%
_{\nabla ,c}$ over $(M,\nabla )$. The metric connection $^{(M)}\widetilde{%
\nabla }$ with respect to ${}\widetilde{g}_{\nabla ,c}$ is as follows:%
\begin{eqnarray*}
^{(M)}\widetilde{\nabla }_{E_{\overline{i}}}E_{\overline{j}} &=&0,\text{ }%
^{(M)}\widetilde{\nabla }_{E_{\overline{i}}}E_{j}=0, \\
^{(M)}\widetilde{\nabla }_{E_{i}}E_{\overline{j}} &=&-\Gamma _{ih}^{j}E_{%
\overline{h}}, \\
^{(M)}\widetilde{\nabla }_{E_{i}}E_{j} &=&\Gamma _{ij}^{h}E_{h}+\frac{1}{2}%
(\nabla _{i}c_{jh}+\nabla _{j}c_{ih}-\nabla _{h}c_{ij})E_{\overline{h}}
\end{eqnarray*}%
with respect to the adapted frame, where $R_{hji}^{\text{ \ \ \ \ }s}$ are
the local coordinate components of the curvature tensor field $R$ of the
symmetric connection $\nabla $ on $M.$
\end{proposition}

\begin{remark}
If ${\nabla _{i}c_{jh}+\nabla _{j}c_{ih}-\nabla _{h}c_{ij}=0}$, the metric
connection $^{(M)}{\widetilde{{}\nabla }}$ in $T^{\ast }M$ of the modified
Riemannian extension $\widetilde{g}_{\nabla ,c}$ coincides with the metric
connection $^{H}\nabla $ of the Riemannian extension $\widetilde{g}_{\nabla }
$, where $^{H}\nabla $ is the horizontal lift of the symmetric connection $%
\nabla $ on $M$.
\end{remark}

For the curvature tensor $^{(M)}{\widetilde{{}R}}$ of the metric connection $%
^{(M)}{\widetilde{{}\nabla }}$, we state the following result.

\begin{proposition}
Let $\nabla $ be a symmetric connection on $M$ and $T^{\ast }M$ be the
cotangent bundle with the modified Riemann extension $\widetilde{g}_{\nabla
,c}$ over $(M,\nabla ).$ The curvature tensor $^{(M)}{\widetilde{{}R}}$ of
the metric connection $^{(M)}{\widetilde{{}\nabla }}$ satisfies the
followings:%
\begin{eqnarray*}
^{(M)}\widetilde{R}(E_{i},E_{j})E_{k} &=&R_{ijk}^{\text{ \ \ \ }h}E_{h}{\,}
\\
&&+\frac{1}{2}\{\nabla _{i}(\nabla _{k}c_{jh}-\nabla _{h}c_{jk})-\nabla
_{j}(\nabla _{k}c_{ih}-\nabla _{h}c_{ik}) \\
&&-R_{ijk}^{\text{ \ \ \ }m}c_{mh}-R_{ijh}^{\text{ \ \ \ }m}c_{km}\}\}E_{%
\overline{h}} \\
^{(M)}\widetilde{R}(E_{i},E_{j})E_{\overline{k}} &=&R_{jih}^{\text{ \ \ \ }%
k}E_{\overline{h}} \\
^{(M)}{\,}\widetilde{R}(E_{i},E_{\overline{j}})E_{k} &=&{\,0,}^{(M)}%
\widetilde{R}(E_{i},E_{\overline{j}})E_{\overline{k}}={\,0,}^{(M)}\widetilde{%
R}(E_{\overline{i}},E_{j})E_{k}=0 \\
^{(M)}\widetilde{R}(E_{\overline{i}},E_{j})E_{\overline{k}} &=&0,^{(M)}%
\widetilde{R}(E_{\overline{i}},E_{\overline{j}})E_{k}=0,^{(M)}\widetilde{R}%
(E_{\overline{i}},E_{\overline{j}})E_{\overline{k}}=0
\end{eqnarray*}%
with respect to the adapted frame.
\end{proposition}

The non-zero component of the contracted curvature tensor field (Ricci
tensor field) $^{(M)}{\widetilde{{}R}}_{\gamma \beta }=$ $^{(M)}{\widetilde{%
{}R}}_{\alpha \beta \gamma }^{\text{ \ \ \ }\alpha }$ of the metric
connection $^{(M)}{\widetilde{{}\nabla }}$ is as follows:%
\begin{equation*}
^{(M)}{\widetilde{{}R}_{jk}=R}_{jk}{.}
\end{equation*}%
For the scalar curvature $^{(M)}\widetilde{r}$ of the metric connection $%
^{(M)}{\widetilde{{}\nabla }}$ with respect to $\widetilde{g}_{\nabla ,c}$ ,
we obtain%
\begin{equation*}
^{(M)}\widetilde{r}=0.
\end{equation*}%
Thus we have the following theorem.

\begin{theorem}
The cotangent bundle $T^{\ast }M$ with the metric connection $^{(M)}{%
\widetilde{{}\nabla }}$ $\ $has a vanishing scalar curvature with respect to
the modified Riemannian extension $\widetilde{g}_{\nabla ,c}$.
\end{theorem}

\bigskip

\end{document}